\documentclass[12pt]{article}

\usepackage{graphicx,array}
\usepackage{amssymb,amsmath}
\usepackage[english]{babel}

\newcommand\ba{\begin{eqnarray}}
\newcommand\ea{\end{eqnarray}}
\newcommand\be{\begin{equation}}
\newcommand\ee{\end{equation}}

\begin{document}

\begin{titlepage}
%\begin{flushright}
%CERN-TH-PH/2012-075
%\end{flushright}
\vskip 1in

\begin{center}

{\Large {\bf From Euler's play with infinite series to the anomalous magnetic moment\footnote{Updated version of Bures-sur-Yvette preprint IHES/P/18/01;   in:	\textit{Quantum Theory and Symmetries with Lie Theory and Its Applications in Physics} Volume 1,  Proceedings of the 10-th International Symposium "Quantum Theory and Symmetries" (QTS10) with 12-th International Workshop "Lie Theory and Its Applications in Physics" (LT12), Varna, Bulgaria, June 2017.}}}

\vspace{4mm}

Ivan Todorov

Institut des Hautes Etudes Scientifiques\\
F-91440 Bures-sur-Yvette, France

and

Institute for Nuclear Research and Nuclear Energy\\
Tsarigradsko Chaussee 72, BG-1784 Sofia, Bulgaria\\
(permanent address) 
e-mail: ivbortodorov@gmail.com

\end{center}

\begin{abstract}

During a first St. Petersburg period Leonhard Euler, in his early twenties, became interested in the Basel problem: summing the series of inverse squares In the words of Andr\'e Weil [W] "as with most questions that ever attracted his attention, he never abandoned it". Euler introduced on the way the alternating "phi-series", the better converging companion of the zeta function, the first example of a polylogarithm at a root of unity.  He realized - empirically! - that odd zeta values appear to be new (transcendental?) numbers.

It is amazing to see how, a quarter of a millennium later, the numbers Euler played with, "however repugnant" this game might have seemed to his contemporary lovers of the "higher kind of calculus", reappeared in the analytic calculation of g-2 - the anomalous magnetic moment of the electron, the most precisely computed and measured physical quantity [K]. 

Mathematicians, inspired by ideas of Grothendieck, are reviving the dream of \'Evariste Galois of
uncovering a group structure in the \textit{ring of periods}, that includes the multiple zeta values and appears in a variety of subjects - from algebraic geometry to Feynman amplitudes. 
 
%Feynman amplitudes are being expressed in terms of a well structured family of special functions and a denumerable set of %numbers - {\it periods}, studied by algebraic geometers and number theorists. The periods appear as residues of the poles of %regularized  primitively divergent multidimensional integrals. In low orders of perturbation theory (up to six loops in the %massless $\varphi^4$ theory) the family of hyperlogarithms and multiple zeta values (MZVs) serves the job. The hyperlogarithms %form a double shuffle differential graded Hopf algebra. 

\end{abstract}

\vfill\eject

\tableofcontents

\end{titlepage}

\begin{center}

\textbf{Introduction}

\end{center}

"The usefulness of useless knowledge" was the provocative title of a 1939 essay of Abraham Flexner, the founding director of the Institute for Advanced Study in Princeton, recalled - and vindicated - in 2017 by the Institute's current director Robert Dijkgraaf. The present story is another illustration of the value of curiosity driven research and of long term thinking in a time full of short term distractions. What could have seemed in mid eighteenth century as Euler's idle play with numbers given by infinite series became a central topic in mathematics and quantum field theory in the last decades.

The paper consists of three distinct parts. The first (Sect. 1) begins with Euler's early encounter of Mangoli's "Basel problem" and follows his repeated assaults on the zeta series and their alternating cousins\footnote{A first hand review of Euler's work and an elementary introduction to multple zeta values is contained in Sects. 1-2 of Cartier's Bourbaki lecture \cite{C01}.}. 
The second (Sect. 2) is concerned with what came to be called the $g-2$ saga (from 1947 to 2017). Some related recent developments in number theory, algebraic geometry and perturbative quantum field theory are surveyed in Sects. 3 and 4.

It should be noted that the stories in Sections 2 and 3 are both open ended. The beginning of a theory of Feynman periods concerns the primitively divergent graphs in a (massless) scalar ($\phi^4$-)theory. It has been only recently realized that $g-2$ calculations (in which infrared divergences are only resolved for gauge invariant sums of Feynman amplitudes) display properties under \textit{Galois coaction}, similar to those discovered in the $\phi^4$ theory. If it should become clear 
from our expos\'e that the study of the Galois coaction on quantum periods is in its infancy, we are not even alluding to the physical problem with the anomalous magnetic moment of the muon whose understanding does not match the success story of g-2 for the electron \cite{J}.  
      
\smallskip
      
\section{Leonhard Euler (1707-1783): zeta values and their alternating companions}
\subsection{The Basel problem}
Having obtained a kind of master degree in Basel at the age of 17, Euler, not quite twenty, got an offer to join the Saint Petersburg Academy. There he was commissioned (fortunately, not for long) into the Russian navy (as he had won, just before leaving Basel, a prize for an essay on ship-building, never having seen a sea-going ship \cite{We}, Chapt. 3, Sect. II). Along with his major work on mechanics (in two volumes), on music theory and naval architecture Euler wrote during his first Petersburg period (1727-41) some 70 memoirs on a great variety of topics.\footnote{A definitive collection of Euler's works, \textit{Opera Omnia}, has been published since 1911 by the Euler Commission of the Swiss Academy of Sciences. By the time of the appearance of his first full scale biography \cite{C15}, at the end of 2015 the edition is nearing completion with over 80 large volumes published. The Enestr\"om index of Euler's papers counts 866 entries.A concise (30-page) biography of Euler with color illustrations is contained in \cite{G}; shorter biographical sketches can be found in \cite{A74, We}.}

 It was early in this period, around 1729, that Euler became first interested in the \textit{Basel problem} - the problem of finding the sum of inverse squares or what we would now call (after Riemann) $\zeta(2)$. It was the dawn of infinite series. One knew how to sum a similar series of inverse rectangles which allowed to prove that $\zeta(2)$ is a real number between one and two:
\ba
\label{zeta2} 
\zeta(s) &=& \sum_1^\infty \frac{1}{n^s}, \, \, \, Re(s)>1,    \nonumber \\
\sum_1^\infty \frac{1}{n(n+1)}=\sum_1^\infty(\frac{1}{n}-\frac{1}{n+1})&=&1<\zeta(2)<1+\sum_1^\infty \frac{1}{n(n+1)} = 2.
\nonumber\\
\ea 
But what exactly is $\zeta(2)$? (It is actually a problem which the young Pietro Mengoli (1625-1686), successor of Cavalieri in Bologna, posed in 1644 \cite{A74} and which excited the brothers-rivals Jacob and Johann Bernoulli in Basel.) Euler first tried to obtain a good numerical estimate for $\zeta(2)$. The series (\ref{zeta2}) is slowly convergent. To get from it a six digit accuracy (by 1731 Euler had $\zeta(2)\simeq 1.644934$ \cite{A}), one would have needed a million terms. On the way of obtaining a faster converging expression Euler first introduced the alternating \textit{phi-series},  
\begin{equation}
\label{phi}
\phi(s) = \sum_{k=1}^\infty (-1)^{k-1}\frac{1}{k^s} = (1-2^{1-s})\zeta(s), \, \quad  \phi(1)=\ln(2),
\end{equation}
a special case of the Dirichlet \textit{L-functions} \cite{S}, introduced over hundred years later. (The series for $\phi(s)$ is convergent for $Re(s)>0$ while the harmonic series $\zeta(1)$ diverges.) More importantly, Euler viewed $\zeta(2)$ as a special value of a power series, the \textit{dilogarithm} \cite{Z}, which also has an integral representation (noted by Leibniz in 1696 in a letter to Euler's teacher, Johann Bernoulli):
\begin{equation}
\label{Li2}
Li_2(x)==\sum_{n=1}^\infty \frac{x^n}{n^2} = -\int_0^x \frac{\ln(1-t)}{t}dt=\int_0^x \int_0^t\frac{dt dt_1}{t(1-t_1)}.
\end{equation}   
$\zeta(2)$ and $\phi(2)$ appear as the values of this function at the two square roots of unity: $\zeta(2)= Li_2(1)$, 
$\phi(2)=Li_2(-1)$. Furthermore, Euler derived the elementary \textit{functional equation} (see \cite{V} Eq. (3)):
\be
Li_2(x)+Li_2(1-x) +ln(x)ln(1-x)=Li_2(1). \,  \, \nonumber
\ee
Setting in it $x=1/2$ he obtained a much faster converging series for $\zeta(2)$. In 1734 Euler had the bold idea to extrapolate Newton's formula for the coefficient of a polynomial in terms of its roots to infinite series (after having obtained an approximation for 
$\zeta(2)$ with twenty decimal places). He announced the beautiful result, $\zeta(2) = \pi^2/6$ in letters to friends (including Daniel Bernoulli) and did not hide his excitement in the introduction to the article (quoted in \cite{V} Sect. 2).  

Here is a sketch of this ingenious direct calculation (see for more detail \cite{B13, V}). Taking the logarithmic derivative of the infinite product expansion of the sine function one finds:
\ba
\label{cot}
cot(x)=\sum_{-\infty}^\infty \frac{1}{x-k\pi} \,  \, \Rightarrow   \nonumber \\
xcot(x) = 1 - 2\sum_{k=1}^\infty\frac{x^2}{k^2\pi^2-x^2} = 1 - 2\sum_{n=1}^\infty \zeta(2n)\left(\frac{x^2}{\pi^2}\right)^n.
\end{eqnarray}
Compared with the power series expansion of $\cot(x)$ this allows to compute $\zeta(2n)$ (at least for small n). In a next assault on the problem Euler recognized the appearance of the numbers $B_n$ encountered in the posthumous work of Jacob Bernoulli (1655-1705) which he called \textit{Bernoulli numbers} and arrived (in 1740) at the beautiful general formula:
\be
\label{zeta2n}
\zeta(2n)=-\frac{B_{2n}}{2(2n)!}(2\pi i)^{2n}, \, \, B_2=\frac{1}{6}, \, \, B_4=-\frac{1}{30}, \, B_6 =
\frac{1}{42}, \,(-1)^{n-1} B_{2n}\in\mathbb{Q}_{>0},
\ee

Euler writes $\zeta(n)=N_n\pi^n$ noting that for even $n\; N_n$ is rational; he has computed $\zeta(3)$ to ten significant figures
and convinced himself that $N_3$ is not a rational number with a small denominator (see \cite{D12} where the original Euler's paper - in Latin - is cited). He conjectured that $N_n$ for odd n might be a function of $ln(2)(=\phi(1))$, \cite{A74}, but this did not work either. Later, in 1749, in his Berlin period, he acknowledges: "for n odd all my efforts have been useless until now" (cited in \cite{A74}). 
In fact, this failure made mathematicians, and especially mathematical physicists, believe that Euler had discovered new transcendental numbers, designed to play an important role in both pure mathematics and quantum field theory. 

\subsection{Memorable mathematical developments}
Trying to find polynomial relations among zeta values Euler was led by the {\it stuffle product}\footnote{A survey of the stuffle and shuffle products and their application to extending the notion of MZV to divergent series may be found in \cite{T, T16}.}
\be
\label{stuf}
\zeta (m) \, \zeta (n) = \zeta (m,n) + \zeta (n,m) + \zeta (n+m)
\ee
to the concept of multiple zeta values (MZVs):
\be
\label{MZV}
\zeta (n_1 , \ldots , n_d) = \sum_{0 < k_1 < \ldots < k_d} \frac{1}{k_1^{n_1} \ldots k_d^{n_d}} .
\ee
which also extends to the alternating phi function (\ref{phi}); for instance,
\be
\label{multiphi}
\phi (m,n) = \sum_{0 < k < \ell}  \frac{(-1)^{k+\ell}}{k^m \ell^n} < 0.
\ee 
In 1731 the 24-year-old Euler introduced the "Euler-Mascheroni constant" (see \cite{La}):
\be
\label{gamma}
\gamma = \lim_{n\rightarrow \infty} (\sum_{k=1}^n \frac{1}{k} - \ln n) =
\sum_{n=2}^\infty (-1)^n \frac{\zeta(n)}{n} (=0.5772...).
\ee

He then discovered the \textit{Euler-MacLaurin formula} which allowed him to express sums in terms of integrals and series of Bernoulli numbers; for instance:
\begin{eqnarray}
\label{EulMcL}
\gamma &=& \sum_{k=1}^n\frac{1}{k} -\int_1^n\frac{dx}{x} -\frac{1}{2n}+\sum_{k=1}^\infty \frac{B_{2k}}{2kn^{2k}}\nonumber\\
&=&
\sum_{k=1}^n\frac{1}{k} -\ln n -\frac{1}{2n} +\frac{1}{12n^2}-\frac{1}{120n^4}+..., \nonumber \\
\zeta(2) &=& \sum_{k=1}^n\frac{1}{k^2} +\int_n^\infty\frac{dx}{x^2} -\frac{1}{2n^2}+\sum_{k=1}^\infty\frac{B_{2k}}{n^{2k+1}}\nonumber\\&=&
\sum_{k=1}^n \frac{1}{k^2} +\frac{1}{n}-\frac{1}{2n^2}+\frac{1}{6n^3}-\frac{1}{30n^5}+ ... .
\end{eqnarray}

Peeling off consecutive prime factors from $\zeta(s)$, starting with two (see \cite{G}),
\be
(1-2^{-s})\zeta(s) = 1 + 3^{-s} + 5^{-s} + ...\,, \nonumber
\ee
Euler discovered in 1737 the fabulous product formula,
\be
\label{primes}
\prod_p (1-p^{-s}) \zeta(s) = 1\ .
\ee

 It was during the subsequent {\it Berlin period} (1741-1766), invited by Frederick II, that, having played with some divergent series, Euler conjectured, in 1749, the functional equation for the $phi$-function, writing (\cite{A74}) ``I shall hazard the following conjecture:
\be
\label{FuncEq}
\frac{\phi (1-s)}{\phi (s)} = - \frac{\Gamma (s) (2^s - 1) \cos \frac{\pi s}2}{(2^{s-1} - 1) \, \pi^s}
\ee
is true for all $s$.'' From  here and from (\ref{phi}) follows immediately the functional equation for $\zeta (s)$, that became, 110 years later, the basis of Riemann's great 1859 paper, \cite{W}.
Euler's work on number theory was done, as Fermat's a century earlier, against a background of contempt towards the field by the majority of mathematicians. He was not deterred. As he once observed "one may see how closely and wonderfully infinitesimal analysis is related not only to ordinary analysis but even to the theory of numbers, however repugnant the latter may seem to that higher kind of calculus" (see \cite{We} Chapt. 3, Sect. V). Euler was repeatedly telling his readers that he felt no need to apologize for the time and effort thus spent, that truth is one, and no aspect of it may be neglected without damage to the whole (\cite{We} Sect. 3.III).

\bigskip

After Euler's death multiple zeta values were all but forgotten for over 200 years before resurfacing simultaneously in quantum field theory and in pure mathematics.

\bigskip
 
\section{The saga of $g-2$}
\subsection{QED and the intrinsic magnetic moment of the electron (1947-1949)}
Divergences in  perturbative quantum electrodynamics (QED) made the founding fathers (Heisenberg, Dirac, Pauli, ...) skeptical about the theory they have created. The younger generation, especially in the US, coming back to the laboratories from a war oriented research, was more pragmatic. Experimentalists were using improved techniques to measure with precision tiny effects in atomic spectra, magnetic resonances and the like; theorists were applying the ideas of Kramers and Pauli-Fierz (1938) on mass renormalization and suggestions of Dirac, Heisenberg and Weisskopf on charge renormalization (also of the 1930's) to explain the new findings by extracting finite corrections from divergent integrals without changing fundamentally the theory (\cite{Sc} - see, in particular, Sect. 4.4, pp. 202-204).  

The magnetic moment of an electron of mass $m$, charge $e$ and spin $\mathbf{s}$  is given by
% \be
% \textbf{\mu} = -g\mu_B\textbf{s}, \;\; \mu_B \,\mathrm{\;being\; the\; Bohr\; magneton} \, \, \mu_B = \frac{e\hbar}{2mc}.
%\ee  
\be
\label{muB}
\mathbf{\mu} = -g\mu_B\textbf{s}, \;\; \mu_B \,\mathrm{\;being\; the\; Bohr\; magneton} \, \, \mu_B = \frac{e\hbar}{2mc}.
\ee  
The constant $g(=g_e)$, called (electron) \textit{g-factor} (or gyromagnetic ratio), is predicted from the Dirac equation to have the value $g=2$.

{\footnotesize \textbf{Aside}: \textit{The role of Gregory Breit}(1899-1981). Back in 1940 Millman and Kusch measure the ratio of the nucleon and electron resonance frequencies in the same magnetic field which is proportional to the corresponding ratio of g-factors. Assuming that the electron magnetic moment takes its Dirac value, $g_e=2$ - the dogma of the day - the authors obtained 
some $0.1\%$ discrepancy in the value of the nuclear magnetic moment compared with other methods. While I.I. Rabi, who has been awarded the Nobel Prize in Physics in 1944 "for his resonance method for recording the magnetic properties of atomic nuclei", and his group at Columbia University were trying to explain the discrepancy by the properties of nuclei, it was Breit who suggested (in a letter to Rabi containing a preview of his letter to the editors of Physical Review) that the origin of the discrepancy may be due to an "intrinsic magnetic moment of the electron" of the order of $\alpha\mu_B$ (where $\alpha$ is the fine structure constant\footnote{introduced by Sommerfeld (1916): $4\pi\epsilon_0\hbar c\alpha= e^2$; in modern particle physics texts the vacuum permittivity $\epsilon_0$ is taken as unity.}). Induced by Rabi, Kusch and Foley found by measuring the Zeeman splitting in Gallium that $g_e$ indeed differed by some two thousands from two, a work which brought the Nobel Prize to Kusch in 1955. Unfortunately, Breit's contribution was not properly acknowledged as witnessed by his subsequent letter to Rabi: "The thing that worries me is ... the general change from the practice of friendliness among scientists. ... I have sent you the letter to the editor, to which Foley and Kusch refer now as follows: 'These results are not in agreement with the recent suggestion made by Breit as to the magnitude of the intrinsic moment of the electron'." (quoted in \cite{Sc} Sect. 5.4, p. 222).}  

Thus,  in 1947, the anomalous magnetic moment of the electron, defined as $a_e:= (g-2)/2$ was found experimentally to be (\cite{K} Sect. 2.1):
\be
\label{aexp47}
               a_e(exp:1947)  = 1.159 (5)\times 10^{-3}.
\ee   
By the end of the year Schwinger computed, using his renormalized QED, the 
one loop electron vertex function in an external magnetic field finding the 
simple result:
\be
\label{ath47}
            a_e(th:1947) = \frac{\alpha}{2\pi}= 1.161... \times 10^{-3}.
\ee

\medskip

$$
\includegraphics[width=5cm]{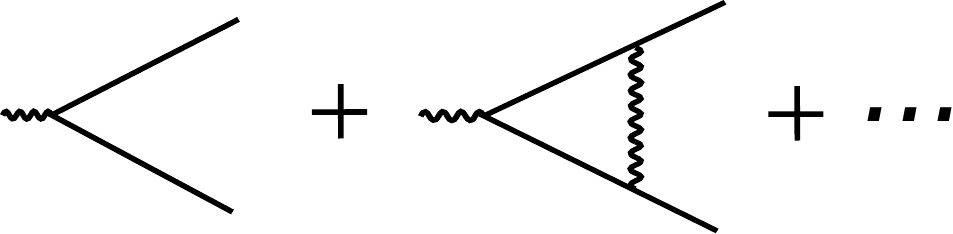}
$$
\label{f1}
\centerline {Figure 1} 
\begin{quote}
The first two terms in the expansion of the electron-photon vertex function.
\end{quote}

\medskip

Around 1949 Feynman had introduced his diagrams and Dyson made them clearer to the community, establishing on the way the equivalence of the approaches of Tomonaga, Schwinger and Feynman. In Feynman's language, Schwingers' calculation amounted to computing a single triangular graph (see Fig.~1). 
%\begin{figure}
%\centering
%\includegraphics[scale=.5]{t1a.pdf}
%\label{f1}
%\caption{}
%\end{figure}
An important simplification came from a joint work of Wheeler (Feynman's PhD adviser) and Feynman who introduced the half sum of the retarded and advanced functions, rediscovering Stueckelberg's \textit{causal propagator}. This makes a calculation that follows Feynman's rules simpler by a factor of two for each line of a diagram (thus eightfold simpler for the diagram on Fig. 1). This simplification opened the possibility for computing higher order corrections.

\subsection{Fourth and sixth order corrections (1949-1996)} 
%\begin{figure}
%\centering
%\includegraphics[scale=.6]{t1b.pdf}
%\label{f2}
%\caption{}
%\end{figure}

\bigskip

$$
\includegraphics[width=12cm]{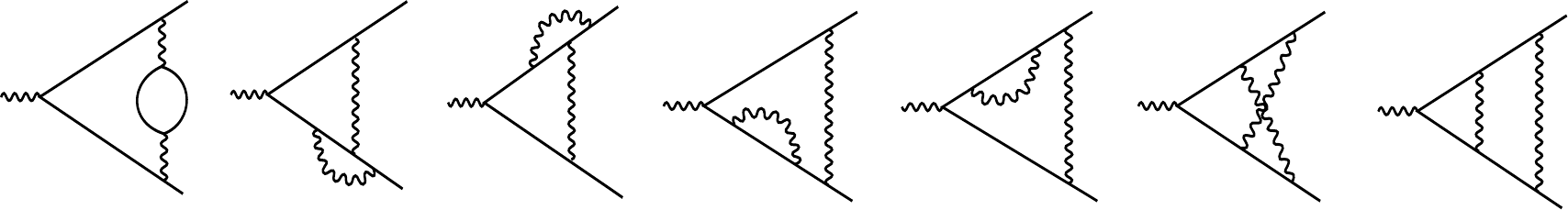}
$$
\label{f2}
\centerline {Figure 2} 
\begin{quote}
Two-loop contributions to the anomalous magnetic moment.
\end{quote}

\medskip

The two-loop calculation of $a_e$ involving the seven diagrams of Fig.~2, a real test of Dyson's renormalization rules, was performed by two young theorists, Karplus and Kroll using the computer at the Institute for Advanced Study in Princeton in 1949 (published in the Physical Review in 1950). It was an influential paper, demonstrating that the work with the newly cooked rules was indeed feasible. When, however, a more precise measurement (by Franken and Liebe) was provided,
\be
\label{exp56}
a_e(exp:1956) =1.1681 (5)\times 10^{-3},
\ee  
it disagreed with the calculation. 

This stimulated Andr\'e Petermann in Geneva to look more closely into the 
six-year-old calculation and he found an error in one of the 
integrals.\footnote{Interviewed some 37 years later (in 1986) Norman Kroll said: 
"[The errors] were arithmetic...  The thing that I learned from that is: in 
doing a complicated calculation,  you have to take the same kinds of precautions 
that an experimenter takes to see that dirt doesn't get in his apparatus.  We 
had some internal checks
but not nearly enough." -\cite{St}} Quite independently Schwinger at Harvard 
found a variational bound which was  violated by the result of Karplus and 
Kroll and assigned the problem of calculating the $\alpha^2$-contribution to 
$a_e$ to his PhD student Charles Sommerfield.\footnote{As Schwinger told at 
the time T.T. Wu. I thank Tai Tsun Wu for this information (contained in his 
email of February 19, 2018).} So Peterman and Sommerfield performed 
separately an \textit{analytic calculation} \cite{Pe, So} and arrived at the same (correct!) 
answer. The result turned out to be nicely expressed in terms of Euler's $\phi$-function (cf. \cite{Sch}):
\be
\label{th57}
a_e(th:1957) = \frac{1}{2}\frac{\alpha}{\pi}+(\phi(3)-6\phi(1)\phi(2) +\frac{197}{2^4 3^2})\left(\frac{\alpha}{\pi}\right)^2 =
1.159 638 (4)\times 10^{-3} .
\ee
Noteworthy, the same weight three combination, $\phi(3)-6\phi(1)\phi(2)$, appears in the second order of the Lamb shift calculation (see \cite{LPR}). 

Comparing the figures and error bars in (\ref{exp56}) and (\ref{th57}) we see that the theoretical predictions have gone ahead of the experimental accuracy. In fact, measurement of g in atomic physics has become stalled. Happily, an entirely different
approach, measurement of $g - 2$ by spin precession in a magnetic field, pursued since ∼ 1953 by the University of Michigan group, has been making a steady progress. After almost 20 years, this method reached the precision of $∼ 3\times 10^{- 6}$ (\cite{K} Sect. 2.4):
\be
\label{exp71}
\, a_e (exp:1971) = 1.159 6577 (35) × 10^{-3} .    \,
\ee  
This is 1400 times more precise than the atomic physics result $a_e(exp:1956)$, forcing theorists to evaluate the three loop (sixth-order) term. That required computing 72 diagrams!

Let me here quote an observer of the "tennis match between theory and experiment" \cite{H}: "Measuring a property of matter with such an extraordinary precision is a labor of years; a single experiment could well occupy the better part of a scientific career. It's not always appreciated that theoretical calculations at this level of accuracy are also arduous and career consuming.... When work on the integrals [corresponding to the 72 Feynman diagrams] got under way in the 1960's it became ... a major impetus to the development of computer algebra systems ... Despite such computational power tools some of the three-loop diagrams resisted analytic result for thirty years. ... It was not until 1995 that a reliable, high precision value of the three-loop contribution was published by Toichiro Kinoshita of Cornell University. He evaluated the 72 diagrams numerically, comparing and combining his results with analytic values that were then known for 67 of the diagrams. A year later the last few diagrams were calculated analytically by Stefano Laporta and Ettore Remiddi of the University of Bologna" \cite{LR}.

 Amazingly, it was again expressed in terms of Euler's multiple $\phi$-values (cf. Eq. (1.2) of \cite{Sch}):
\begin{eqnarray}
\label{th96}
a_e(th:1996) &= &\frac12 \, \frac{\alpha}{\pi} + \left[ \phi (3) - 6 \, \phi (1) \, \phi (2) + \phi (2) + \frac{197}{2^4 \, 3^2} \right] \left( \frac{\alpha}{\pi} \right)^2 \nonumber \\
+& &\!\!\!\!\!\!\!\Biggl[ \frac2{3^2} \, (83 \, \phi (2) \, \phi (3) - 43 \, \phi (5)) - \frac{50}3 \, \phi (1,3) + \frac{13}5 \, \phi (2)^2  \\
+& &\!\!\!\!\!\!\!\frac{278}3 \left( \frac{\phi (3)}{3^2} - 12 \, \phi \, (1) \, \phi (2) \right) + \frac{34202}{3^3 \, 5} \, \phi (2) + \frac{28259}{2^5 \, 3^4} \Biggl] \left( \frac{\alpha}{\pi} \right)^3 + \ldots \nonumber\\
=&& \!\!\!\!\!\!\! 1.159652201(27)\times10^{- 3} \nonumber
\end{eqnarray} 
(the error reflecting the uncertainty in $\alpha$). This beats the accuracy of the Michigan experiment. But in the meantime experimenters have not stayed idle. While the spin precession method hits the ceiling, an approach that utilizes the spin
and cyclotron resonances in a Penning trap (which began around 1958) was being pursued by the group of Dehmelt\footnote{Hans Dehmelt (1922-2017) shared the 1989 Nobel Prize in Physics "for the development of the ion trap technique".} et al. at the University of Washington.

After 30 years this approach led to three orders of magnitude improvement over the precession measurement of the Michigan group. Their results for an electron and a positron are:
\ba
\label{exp87}
a_{e^-}(exp :1987) &=& 1.159 652 1884 (43) \times 10^{-3}, \nonumber\\
a_{e^+}(exp :1987) &=& 1.159 652 1879 (43) \times 10^{-3}.
\ea
It is hard to overestimate the beauty and the significance of a formula like (\ref{th96}) given the precision with which it is confirmed experimentally. Dyson \cite{D52} has given an argument that the perturbative expansion in QED is likely to be divergent. Its close agreement with experiment, on the other hand, makes us believe that it is asymptotic. Individual terms have then a meaning of their own, both as special exactly known numbers and as measured quantities. Subsequent work, to be reviewed in Sect. 2.3 below, confirm the expectation that higher powers of $\frac{\alpha}{\pi}$ provide at least a hundred times smaller contribution. One is tempted to place these formulas among what Salviati (the \textit{alter ego} of Galileo) elevates to "those few which the human intellect does understand, I believe that its knowledge equals the Divine in objective certainty" \cite{Ga}

\subsection{Eight's order and beyond (1996-2017)}
The result (\ref{exp87}) means that theory must be extended to four-loop graphs (or to the eighth-order in e) since 
$(\frac{\alpha}{\pi})^4 \sim 29\times 10^{-12}$. This required  computing 891 diagrams. A thousandfold increase of the computer time has brought a thirty-fold improvement in precision. After more than twenty years of hard work Kinoshita and coworkers achieved in 2012 a reliable numerical estimate for the coefficient to $(\frac{\alpha}{\pi})^4$ (see \cite{K} and references therein\footnote{Kinoshita uses $A_1^{(2n)}$ instead of $a_n$.}):
\be
\label{a4num}
a_e = \sum_{n>0} a_n (\frac{\alpha}{\pi})^n, \, \, \, a_4 = −1.9106 (20) .
\ee

In the words of the spectator \cite{H} "Attacking all those intricately tangled diagrams by analytic methods is hopeless for now." Yet, it was achieved\footnote{Is it possible that the hero of these (and other) calculations, Stefano Laporta, never had a tenure in Bologna (as I learned from David Broadhurst, December, 2014)?} \cite{L17}!
The result involves hyperlogarithms at twelfth roots of unity (called multiple Deligne values in \cite{B14}), one-dimensional integrals of products of complete
elliptic integrals and six finite parts of master integrals, evaluated up to 4800 digits. The outcome,
\be
\label{aLap}
a_4 = -1.912245764926445574152647167439830054060873390658725...
\ee
recovers the numerical computation (\ref{a4num}) within its uncertainty. In the meantime, in 2008, the Harvard group of Gabrielse et al. developing the cylindrical Penning trap method increased 15 times the accuracy of the University of Washington group and obtained the best now available measurement:
\be
\label{exp08}
a_e(exp:2008) = 1.159 652 180 73 (28) \times 10^{-3}.
\ee 

With such a precision, one must also take into account heavier particles such as muon, tau, hadrons and weak bosons in the intermediate states. This has been done both numerically and analytically (as reviewed in \cite{K}). The group of the old veteran Toichiro Kinoshita (born in January, 1925) also estimated the contribution of the 12 672 five-loop diagrams. His result of 2012 agrees with the refined version of Laporta (2017) within the respective errors:
\ba
\label{thKL}
a_e(th:2012) &=& 1.159 652 181 78 (6)(4)(2)(77) \times 10^{-3}, \nonumber\\
 a_e(th:2017) &=& 1.159 652 181 664(23)(16)(763)\times 10^{-3}.
\ea
The uncertainties in the first result (of 2012) come from 4-loop, 5-loop, hadronic and weak contributions and  - the largest one - from the measurement of the fine structure constant $\alpha$; in the second figure there is no 4-loop uncertainty and again the largest error comes from $\alpha$. Kinoshita estimates the discrepancy between theory and experiment as
\begin{equation}
\label{diff}
\, a_e(exp:2008) \, - \, a_e(th:2012) = −1.05 (82) \times 10^{-12} .
\end{equation}

It becomes even smaller with the Laporta result. As the biggest uncertainty comes from the fine structure constant, it makes sense, conversely, to determine $\alpha$ from the measurement of $a_e$ and compare it with the best other value (determined by measuring the recoil velocity of a rubidium atom when it absorbs a photon):
\ba
\label{alpha}
\alpha^{-1}(a_e ) &=& 137.035 999 1727 (68)(46)(19)(331), \nonumber\\
 \, \alpha^{-1}(Rb11) &=& 137.035 999 049 (90).
\ea  
The biggest uncertainty (331) in $\alpha^{− 1}(a_e )$ comes from the measurement of $a_e(exp:2008)$. Laporta replaces the last three figures (727) by 596; the first uncertainty (66), corresponding to four loops, is absent in his case, while the next two are substituted by (27) (for the five loop term) and (18) (for the hadronic and electroweak corrections).

The unprecedented agreement between theory and experiment provides a more solid confirmation of the renormalization procedure 
in QED than was expected by its creators - as witnessed by the following letter by Freeman Dyson to Gerald Gabrielse:
"... As one of the inventors, I remember that we thought of QED in 1949 as a temporary and a jerry-built structure, with mathematical inconsistencies and renormalized infinities swept under the rug. We did not expect it to last more than 10 years
before some more solidly built theory would replace it... Now 57 years have gone by and that ramshackle structure still stands... It is amazing that you can measure her dance to one part per trillion and find her still following our beat..." (cited in \cite{GH}).

\section{Periods in quantum field theory}
Euler played with rather special numbers - which since have names: (multi)zeta and phi values, powers of $\pi, \ln n$, ... . Why should we care? If you take "at random" a number on the real axis, it is likely to be a nameless transcendental. It turns out that numbers in quantum field theory (QFT) - the basis of particle physics - are not random. They are \textit{periods} \cite{BB, BW, B15} - a countable set of complex numbers singled out by mathematicians \cite{KZ, A08, M-S}.   
%The functional $res\, G$ is a {\it period} according to the definition of \cite{KZ, M-S}.

This is a place to pose and look at some developments in number theory, algebraic geometry and in perturbative quantum field theory during the last decades which indicate that the overlap between these subjects is not fortuitous.

\subsection{The ring of periods}
The multiple zeta values and their generalizations are special values of hyperlogarithms which find their natural playground among Chen's iterated path integrals \cite{C}. Their systematic study in mathematics was resurrected in the work of Zagier \cite{Z92} and Hoffman \cite{H92} and shortly after in particle physics by Broadhurst et al. (see \cite{BK} and references to earlier work cited there). The notion of a period crystallized in the work of Kontsevich and Zagier at the turn of the century \cite{KZ}. It looks deceptively simple: a complex number is a period if its real and imaginary parts can be written as absolutely convergent integrals of rational functions in domains given by polynomial inequalities with all coefficients in $\mathbb{Q}$. Every algebraic number is a period. For instance, the n-th root of the positive integer k can be defined by
\begin{equation}
\label{nroot}  
k^{1/n} =  \int_{0<x, x^n<k}dx.
\end{equation}
 Moreover, the set $\mathbb{P}$ of all periods would not change if we replace everywhere in the definition "rational" by "algebraic". If we denote by $\bar{\mathbb{Q}}$ the {\it field of algebraic numbers} (the inverse of an algebraic number being also algebraic) then we would have the inclusions
\begin{equation}
\label{QP}
\mathbb{Q}\subset \bar{\mathbb{Q}} \subset \mathbb{P} \subset \mathbb{C}.
\end{equation}
The periods form a ring (they can be added and multiplied) but the inverse of a period needs not be a period. The set of all periods is still countable although it contains infinitely many transcendental numbers, including
\begin{equation}
\label{Periods}
\pi = \iint_{x^2+y^2\leq 1} dx dy, \, \, \ln n = \int_1^n \frac{dx}{x}, \, \, n=2, 3, ...,
\end{equation}
as well as the values of iterated integrals, studied in \cite{C, B09}, at algebraic arguments. They include the classical multiple zeta values and their generalization, the multiple Deligne values, mentioned in Sect. 2.3. The basis $e$ of natural logarithms, the Euler constant $\gamma = -\Gamma'(1)$, as well as $\ln(\ln n), \ln (\ln (\ln n))$, ..., and $1/\pi$ are believed (but not proven) not to be periods. For euclidean momenta (renormalized) Feynman amplitudes in an arbitrary (relativistic, local) QFT can be normalized in such a way that for rational values of the coupling constants and ratios of dimensional parameters they are periods. (It is proven in \cite{BB} that Feynman amplitudes with rational parameters  can  be  written as  a  product of a  Gamma factor and a meromorphic function such that the coefficients of its Taylor expansion are all periods; a similar result, using the Laurent expansion of dimensionally regularized amplitudes is obtained in \cite{BW}.) Defining the \textit{graph polynomial} of a connected Feynman graph $G$ by 
\begin{equation}
\label{GraphPoly}
\Psi_G(\alpha) = \sum_{T\subset G} \prod_{e\notin E(T)}\alpha_e,
\end{equation}
where the sum is over the spanning trees $T\subset G$ with a set of edges $E(T)$, one can write the period of a primitively logarithmically divergent graph in a massless theory as the projective integral\footnote{More generally one is dealing with integrands depending on external momenta and masses, again expressed as ratios of polynomials \cite{B15} (explained in Panzer's thesis \cite{P}).}:
\be
\label{eq28}                  
P (G) = \int_{0<\alpha_e<\infty}\frac{\Omega(G)}{\Psi_G(\alpha)^2}, \, \,
\Omega(G) = \sum_{i=1}^n (-1)^{n-i} \alpha_i d\alpha_1...\hat{d\alpha_i} ... d\alpha_n
\ee
where the hat over a factor means that this factor should be omitted in the product. The same numbers appear in the expansion of the renormalization group beta function (see \cite{S97, GGV, S18}). Here $n =
|E(G)|$ is the number of edges in the graph $G$ which, for a logarithmically divergent amplitude, equals twice the number of loops of the graph: $n=2h(G)$ (in general, for a connected graph with $|E|$ edges and $|V|$ vertices $h=|E|-|V|+1$; it gives the \textit{first homology class} (Betti number) of $G$). As it stands, the integral $P(G)$ is logarithmically divergent. It can be defined as the \textit{residue} of the pole with respect to some analytic regularization \cite{NST}. It is equal to the restriction of the integral (\ref{eq28}) to any transverse surface, say, $\alpha_n=1$ (as in \cite{PS}). There is a number of equivalent expressions for the period of a primitive graph (as an integral in momentum, position or Schwinger parameters space - see \cite{Sch} for the case of the massless $\phi^4$ theory). With the above normalization (adopted in \cite{Sch, BrS, PS} among many others) it yields rational residues for one- and two-loop graphs.
For graphs with three or higher number of loops $h$, one encounters, in the $\phi^4$ theory, multiple zeta values of overall weight not exceeding $2h - 3$ (cf. \cite{BK,Sch,S14}). In fact, the only periods at three, four and five loops (in the $\phi^4$ theory) are integer multiples of $\zeta(3), \zeta(5)$ and $\zeta(7)$, respectively. The first double zeta value, $\zeta(3, 5)$, appears at six loops (in the combination $\frac{2}{5}(29\zeta(8)-12\zeta(3,5)) -9\zeta(5)\zeta(3)$ with an integer coefficient - see the census in \cite{Sch} or the graph $P_{3,5}$ of Sect 2.2 of \cite{B15}). All {\it known} residues were (up to 2013) rational linear combinations of multiple zeta values (MZVs) \cite{BK,Sch}. The seven loop graph whose completion (with a point at infinity) is displayed on Fig. 3 was demonstrated in 2014 \cite{P,B14, PS} to involve values of {\it hyperlogarithms} at sixth roots of unity (a special case of multiple Deligne values). The story does not seem to approach an end: mathematicians and physicists are already happily exploring - and realizing in perturbative QFT - elliptic polylogarithms\footnote{related to Gelfand-Kapranov-Jelevinsky hypergeometric series and to mirror symmetry see \cite{V18}} (for different mass scales) \cite{BL, ABW, BDDT} and modular forms (starting with eight loops) \cite{BrS, B13, B17}. Happily, as we shall see, one is also uncovering some general patterns.  

\textit{Remark on the terminology.} Yves Andr\'e \cite{A08} calls the above defined periods \textit{effective periods} and reserves the term period for the elements of the algebra $\mathbb{P}[\frac{1}{2\pi i}]$, denoted by $\hat{\mathbb{P}}$ and called the \textit{algebra of abstract periods} in Sect. 4 of \cite{KZ}. It is noteworthy that Feynman periods do not require inverting $\pi$.  

$$
\includegraphics[width=3cm]{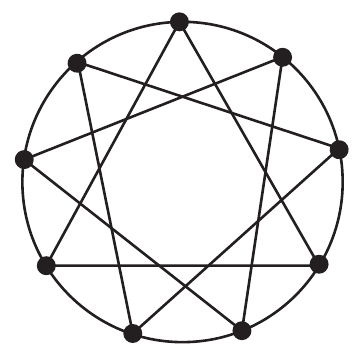}
$$
\label{f3}
\centerline {Figure 3} 
\begin{quote}
The completion (defined in \cite{S14,PS}) of a 7-loop graph whose period involves multiple polylogarithms at sixth roots of unity.
\end{quote}

\subsection{A Galois group for periods?}
The study of the field $\bar{\mathbb{Q}}$ of algebraic numbers, the set of roots of polynomial equations with rational coefficients, "evolved in symbiosis with Galois theory" (to cite \cite{A08}). The question whether this idea can be extended to some class of transcendental numbers had occurred to Galois himself\footnote{"Mes principales m\'editations
depuis quelque temps \'etaient dirig\'ees sur l'application \`a  l'analyse transcendante de la th\'eorie de l'ambiguit\'e." - see \cite{A}.}. The problem of studying the 
Galois group of (motivic) periods was clearly stated in \cite{A08} and vigorously pursued and mutually stimulated in the work of Francis Brown \cite{B11, B12, B15, B15-17}, Oliver Schnetz 
and Erik Panzer \cite{BrS, S14, BS, P, PS, S17}. Our introductory remarks below only aim to attract the readers' interest to this work. Brown \cite{B15} introduces the notion of 
motivic period to an audience of mathematical physicists with the simplest example of $2\pi i= \int_\gamma\frac{dx}{x}$ where $\gamma$ is a small loop which winds in the positive 
direction around the origin in the complex plane. The integral  won't change if we add to the closed form $\frac{dx}{x}$ any exact one-form. It only depends on the \textit{de Rham 
cohomology} class $H^1_{dR}$.  Similarly, the domain of integration is a closed chain and, again by Stokes' theorem, can be modified by the boundary of a two-chain. Hence we only care about its class in the quotient which is the \textit{singular Betti homology} dual to $H^1_B$. Thus a motivic version $(2\pi i)^\mathfrak{m}$ of $2\pi i$ represents a triple
\begin{equation}
\label{Motriple}
(2\pi i)^\mathfrak{m} = \{ [\frac{dx}{x}]\in H^1_{dR}, [\gamma]\in (H^1_B)^\vee; comp\}
\end{equation}  
where the \textit{comparison map}  $comp$ stands for an isomorphism of the complexification of the two cohomology spaces (originally defined algebraically over the rational numbers $\mathbb{Q}$). The question arises: why does one need such an abstract notion as \textit{motivic periods} for a set of real or complex numbers given by convergent integrals of rational functions? To give a "working man answer" to it, we  start with the picture of MZVs as special values of hyperlogarithmic functions which form a differential graded Hopf algebra (see e.g. \cite{T, T16} for a review and references). In particular, the classical polylogarithms are equipped with the coproduct $\Delta^{it}$ of iterated integrals:
\begin{equation}
\label{coprodLn}
\Delta^{it} Li_n(z)= Li_n(z)\otimes 1 +\sum_{k=0}^{n-1}\frac{(\ln z)^k}{k!}\otimes Li_{n-k}(z).
\end{equation}     
Suitably interpreted, this coproduct encodes both the differential equation $dLi_n(z)=(d\ln z)Li_{n-1}(z)$ and the monodromy of the polylogarithm around $z=1$. It is therefore not surprising that such Hopf algebras find many applications in the study of Feynman amplitudes (see e.g. \cite{D} as well as the recent papers \cite{ABDG, BDDT} and references to earlier work cited there). This coproduct, however, does not work as it stands for $z=1, Li_n(1) = \zeta(n)$, since it would not preserve, for instance, the numerical relation $2\zeta(2)^2 = 5\zeta(4)$. Moreover, the hyperlogarithms form a \textit{graded} (double shuffle) algebra while for their specialization, the MZVs, such a property requires a widely open conjecture of linear independence of MZVs of different weights over the rationals (and, furthermore, even algebraic independence of $\zeta(2n+1)$ and $\pi$). To retrieve the above structures we substitute, in particular, $\zeta(2)$ with a \textit{motivic} period,\footnote{Motivic periods were first associated with graph polynomials in \cite{BEK}. An accessible review of formal (and motivic) double zeta values is contained in \cite{C12}.} $\zeta^\mathfrak{m}(2)$, and assume that rather than a coproduct it satisfies a \textit{coaction}\footnote{There are, unfortunately, two opposite conventions of writing this coaction, right and left, in \cite{B15, B17} and in \cite{PS}, respectively. We adopt that of \cite{PS} which coincides with the one of our earlier work \cite{T, T16}.} - and the same is true for 
$(2\pi i)^\mathfrak{m}$:
\begin{equation}
\label{Delta2}
\Delta \zeta^\mathfrak{m}(2) = 1\otimes \zeta^\mathfrak{m}(2), \, \Delta (2\pi i)^\mathfrak{m} = 1\otimes (2\pi i)^\mathfrak{m}. 
\end{equation}    
More generally, for, say, the motivic version of the dilogarithm we write:
\begin{equation}
\label{DeltaDilog}
\Delta Li_2^\mathfrak{m}(z) = 1\otimes Li_2^\mathfrak{m}(z) + \ln^\mathfrak{u}(z)\otimes Li_1^\mathfrak{m}(z) +
Li_2^\mathfrak{u}(z)\otimes 1, 
\end{equation}
where $Li_2^\mathfrak{m}(1) = \zeta^\mathfrak{m}(2)$, while $\ln^\mathfrak{u}(1)=Li_2^\mathfrak{u}(1)=0$ and could be called \textit{unipotent de Rham} periods \cite{B15}. Odd zeta values obey the simple coaction rule 
\begin{equation}
\label{zetaOdd}
\Delta\zeta^\mathfrak{m}(2n+1) = 1\otimes \zeta^\mathfrak{m}(2n+1)+\zeta^\mathfrak{u}(2n+1)\otimes 1 
\end{equation}
(similar to the coproduct  of primitive elements of a Hopf algebra). Thus, the motivic periods $\mathbb{P}^\mathfrak{m}$ appear as a comodule over the Hopf algebra $\mathbb{P}^\mathfrak{u}$ of unipotent de Rham periods.  
Their relation to the iterative integrals and the (real) MZVs is encoded in the existence of a surjective \textit{period map}:
\begin{equation}
\label{per}
per Li_n^\mathfrak{m}(z) = Li_n(z) \Rightarrow  per \zeta^\mathfrak{m}(n) = \zeta(n).
\end{equation} 
The widely open "standard conjectures" would follow if this map were also injective - and hence an isomorphism. The \textit{Galois group for motivic periods} is dual to the Hopf algebra $\mathbb{P}^\mathfrak{u}$.
Thus, the \textit{Galois conjugates} to a motivic period are the right entries in the coproduct of this period. In particular, while the set of conjugates of $\zeta^\mathfrak{m}(2n)$ is one-dimensional (it consists of rational multiples of $\zeta^\mathfrak{m}(2n)$), the $\mathbb{Q}$-space of conjugates of odd zeta values is two-dimensional, spanned by $\zeta^\mathfrak{m}(2n+1)$ and 1.  

As observed in \cite{PS} Feynman periods, in particular $\mathbb{P}_{\phi^4}$ periods (corresponding to the $\phi^4$ theory), are sparse.
For instance, no Feynman graph is known to give rise to the simplest zeta value, $\zeta(2)$. Combined with the highly non-trivial observation in \cite{PS} that $\mathbb{P}_{\phi^4}$ is stable under the Galois coaction, elevated to a conjecture about the motivic $\mathbb{P}_{\phi^4}$ periods, implies that one would never encounter products of the type $\zeta(2)\zeta(2n+1)$ in $\mathbb{P}_{\phi^4}$ (since $\zeta(2n+1)^\mathfrak{u}\neq 0$). What makes such observations even more significant is the result of Brown that motivic Feynman periods of a given \textit{type} (determined by the number of external lines and different masses) span a finite dimensional space (Theorem 5.2 of \cite{B15}). This theorem allows to predict the type of periods of given weight (that generalizes the notion of weight of MZVs) in amplitudes of any order. An illustration of what this means is the observation in \cite{Sch} that the period of a six-loop graph, $P_{3,5}$ (encountered in Sect. 3.1 above), also appears in a seven-loop period (multiplied by $252\zeta(3)$ - see Eq. (6.2) of \cite{B15}). 

Using the amazing result of \cite{L17}, Schnetz  \cite{S17} demonstrates that the Galois coaction structure in QED is similar to the one conjectured for the $\phi^4$ theory - with small dimensions of the 
$\mathbb{Q}$ vector spaces of Galois conjugatres of the $g-2$ periods. In contrast to the way G. Breit has been treated back in 1947 (Sect. 2.1) Schnetz writes: "In this note we process a result by S. Laporta. We keep this note short to emphasize
that the results should be mostly attributed to S. Laporta."

\section{Outlook}
\setcounter{equation}{0}
Analytic methods in perturbative quantum field theory have been developed for over seventy years, starting with the early success of QED in explaining newly observed deviations from the predictions of the Dirac equation, and, in particular, after the first correct two-loop computation of the anomalous magnetic moment of the electron in 1957. After the advent of LHC there have been impressive advances in the calculation of on-shell multileg scattering amplitudes (which can be traced back from \cite{ABHY, C17, FL}) and one loop graphs \cite{AY, BDDT}. Unfortunately, they do not enhance our understanding of the all-important calculation of the anomalous magnetic moment, where the simplicity of the outcome - Eq. (17) and even Eq. (19) - contrasts with the difficulty of its derivation. Should not there be a more direct way to calculate such a basic physical quantity that is a pure number (with no parameters involved!)? 

Closer to the discussion in Sect. 3, various seemingly well founded conjectures concerning the general type of functions and associated numbers that appeared in such calculations had to be revised, however, with each subsequent loop order (or extra mass parameter). As hyperlogarithms and their values at roots of unity were found not sufficient to express all massive and higher order Feynman amplitudes, mathematical physicists started exploring elliptic polylogarithms and modular forms (as surveyed in Sect. 3.1). It was at this point that Brown \cite{B15, B15-17}, stimulated by the tirelessly pursued calculations by Panzer and Schnetz \cite{PS}, pushed forward the visionary idea of Galois, Grothendieck (Cartier \cite{C98}, Connes et al. \cite{CM}, Andr\'e \cite{A08}, ...) of a \textit{cosmic Galois group} in QFT. Remarkably, the Galois coaction principle seems to work not just for contributions of individual Feynman diagrams, but also in the $g-2$ calculation \cite{L17, S17} where an infrared finite result is only obtained for gauge invariant sums of graphs of a given order, and for the "hexagon bootstrap amplitude" \cite{CDMH}  obtained  by combining perturbative input with the axiomatic Steinmann's relations (see also \cite{D17}). 

%The dimensions of weight spaces of MZVs (which exhaust the Feynman periods up to six loops in the massless $\varphi^4$ theory) %do not exceed - and are conjectured to coincide with - their motivic counterparts studied by Francis Brown \cite{B11}. Values %of hyperlogarithms at sixth roots of unity first appear at seven loops. For the two-loop sunrise integral with massive %propagators one encounters multiple elliptic polylogarithms \cite{BV, ABW, BKV}.

The interplay between algebraic geometry, number theory and perturbative QFT, that vindicates Euler's "useless efforts", is a young and vigorous subject and our survey is far from complete. We have not touched upon the application of cluster algebras to multileg on-shell Feynman amplitudes - see \cite{GGSVV}. Nor did we refer to the application of motivic MZVs to string perturbation theory, cited in \cite{B15, B15-17, B17}.     

Brown's concluding remarks in \cite{B15} are inspiring:
The motivic Galois group of hidden symmetries provides an organizing principle for the structure of Feynman amplitudes, 
valid to all orders in perturbation theory. We are only just beginning to scratch the surface of this structure.
 
\smallskip

It is a pleasure to thank Pierre Cartier for his insightful remarks. 

{\footnotesize The author thanks IHES and the National Center for Competence in Research "The Mathematics of Physics" (NCCR SwissMAP) for hospitality during the final stage of this work. He also acknowledges the help of Ludmil Hadjiivanov and Mikhail Stoilov.}

\end{document}